\documentclass[reqno,12pt]{amsart}
\usepackage{amssymb}
\usepackage{amsmath}
\usepackage{amsfonts}

\theoremstyle{plain}

\numberwithin{equation}{section}
\input{tcilatex}

\begin{document}
\title[Some properties of Lipschitz strongly $p$-summing operators]{Some
properties of Lipschitz strongly $p$-summing operators}
\author{Khalil SAADI}
\address{\"{y} }
\date{}
\dedicatory{University of M'sila, Laboratoire d'Analyse Fonctionnelle et G%
\'{e}om\'{e}trie des Espaces, 28000 M'sila, Algeria\\
kh\_saadi@yahoo.fr}

\begin{abstract}
We consider the space of molecules endowed with the transposed version of
the Chevet-Saphar norm and we identify its dual space with the space of
Lipschitz strongly $p$-summing operators. We also extend some old results to
the category of Lipschitz mappings and we give a factorization result of
Lipschitz $\left( p,r,s\right) $-summing operators.
\end{abstract}

\maketitle

\begin{center}
\textsc{0. Introduction}
\end{center}

An area of research that is currently very active is the study of nonlinear
geometry, of Banach spaces or even of general metric spaces, by borrowing
ideas and insights from the linear theory of Banach spaces. A very powerful
tool in the latter is given by the class of $p$-summing operators, so that
naturally has led several authors to investigate Lipschitz versions of them
starting with the seminal paper $\left[ \text{7}\right] $. J.A. Chavez
Dominguez $\left[ \text{2}\right] $ has explored more properties of this
class and has defined a norm on the space of molecules of which dual space
coincides with the space of Lipschitz $p$-summing. The aim of this paper is
to continue to study the same ideas developed in $\left[ \text{2}\right] .$
We try \ to give a similar treatment to the class of Cohen strongly $p$%
-summing operators. We consider the transposed version of the norm of
Chevet-Saphar and we show that the dual of the space of molecules endowed
with this norm coincides with the space of Lipschitz-Cohen strongly $p$%
-summing. Some old results have been established, namely a version of
Grothendieck's theorem and the relationship between the Lipschitz mapping $%
T:X\rightarrow Y$ and its linearization $\widehat{T}:\mathcal{F}\left(
X\right) \longrightarrow Y$ for the concepts of $p$-summing and Cohen
strongly $p$-summing.

The paper is organized as follows.

First, we recall some standard notations which will be used throughout. In
section 1, we define a norm on the space of molecules that is inspired by
the Chevet-Saphar norms. We give and prove an integral characterization of a
linear form on this new space. In section 2, we give the definition of
Lipschitz-Cohen strongly $p$-summing for maps from a metric space to a
Banach space. We show that this space is precisely the dual of the space of
molecules described in section 1. Finally, in section 3 we study some basic
properties of these Lipschitz-Cohen strongly $p$-summing operators, drawing
parallels to the linear theory. Some interesting results have been obtained
namely the Grothendieck theorem and the relationship between the Lipschitz
mapping and its linearization for certain concept of summability. The last
part of this section is devoted to study a factorization result of Lipschitz 
$\left( p,r,s\right) $-summing operators like the one given in linear case.
We show that the map $T$ is Lipschitz $\left( p^{\ast },r,s\right) $-summing
if, and only if, $T$ can be written as $T=T_{2}\circ T_{1}$ where $T_{1}$ is
Lipschitz $r$-summing and $T_{2}$ is Lipschitz-Cohen strongly $s^{\ast }$%
-summing.%
\vspace{0.5cm}%

Now, we recall briefly some basic notations and terminology. Throughout this
paper we always consider metric spaces with a distinguished point (pointed
metric spaces) which we denote by $0$. Let $X$ be a pointed metric space. We
denote by $X^{\#}$ the Banach space of all Lipschitz functions $%
f:X\longrightarrow \mathbb{R}$ which vanish at $0$ under the Lipschitz norm
given by%
\begin{equation*}
Lip\left( f\right) =\sup \left\{ \frac{\left\vert f\left( x\right) -f\left(
y\right) \right\vert }{d\left( x,y\right) }:x,y\in X,x\neq y\right\} .
\end{equation*}%
Consider $1\leq p\leq \infty $, for sequences of the same length $\left(
\lambda _{i}\right) _{i=1}^{n}$ of real numbers and $(x_{i})_{i=1}^{n}$, $%
\left( x_{i}^{\prime }\right) _{i=1}^{n}$ of points in $X$, we denote their
weak Lipschitz $p$-norm by%
\begin{equation*}
w_{p}^{Lip}\left( \left( \lambda _{i},x_{i},x_{i}^{\prime }\right)
_{i}\right) =\sup_{f\in B_{X^{\#}}}(\sum_{i=1}^{n}\left\vert \lambda
_{i}\left( f\left( x_{i}\right) -f\left( x_{i}^{\prime }\right) \right)
\right\vert ^{p})^{\frac{1}{p}}.
\end{equation*}%
We denote by $\mathcal{F}\left( X\right) $ the free Banach space over $X$,
i.e., $\mathcal{F}\left( X\right) $ is the completion of the space%
\begin{equation*}
AE=\left\{ \sum_{i=1}^{n}\lambda _{i}m_{x_{i}x_{i}^{\prime }},\left( \lambda
_{i}\right) _{i=1}^{n}\subset \mathbb{R},\text{ }\left( x_{i}\right)
_{i=1}^{n},\left( x_{i}^{\prime }\right) _{i=1}^{n}\subset X\right\} ,
\end{equation*}%
with the norm%
\begin{equation*}
\left\Vert m\right\Vert _{\mathcal{F}\left( X\right) }=\inf \left\{
\sum_{i=1}^{n}\left\vert a_{i}\right\vert d\left( x_{i},x_{i}^{\prime
}\right) :m=\sum_{i=1}^{n}a_{i}m_{x_{i}x_{i}^{\prime }}\right\} ,
\end{equation*}%
where the function $m_{x_{i}x_{i}^{\prime }}:X^{\#}\rightarrow \mathbb{R}$
is defined as follows%
\begin{equation*}
m_{x_{i}x_{i}^{\prime }}\left( f\right) =f\left( x_{i}\right) -f\left(
x_{i}^{\prime }\right) .%
\vspace{0.5cm}%
\end{equation*}%
We have%
\begin{equation*}
\mathcal{F}\left( X\right) ^{\ast }=X^{\#}.
\end{equation*}%
For the general theory of free Banach spaces, see $\left[ \text{1, 8, 9, 12,
16}\right] $. Let $X$ be a pointed metric space and $Y$ be a Banach space,
we denote by $Lip_{0}\left( X;Y\right) $ the Banach space of all Lipschitz
functions $T:X\rightarrow Y$ such that $T\left( 0\right) =0$ with pointwise
addition and Lipschitz norm. We also denote by $\mathcal{F}\left( X;Y\right) 
$ the vector space of all $Y$-valued molecules on $X,$ i.e.,%
\begin{equation*}
\mathcal{F}\left( X;Y\right) =\left\{ \sum_{i=1}^{n}\lambda
_{i}y_{i}m_{x_{i}x_{i}^{\prime }},\left( \lambda _{i}\right) _{i}\subset 
\mathbb{R},\left( y_{i}\right) _{i}\subset Y,\text{ }\left( x_{i}\right)
_{i},\left( x_{i}^{\prime }\right) _{i}\subset X\right\} .
\end{equation*}%
For any $T\in Lip_{0}\left( X;Y^{\ast }\right) ,$ we denote by $\varphi _{T}$
its correspondent linear function on $\mathcal{F}\left( X;Y\right) $ defined
by%
\begin{equation*}
\left\langle \varphi _{T},m\right\rangle =\left\langle T,m\right\rangle ,
\end{equation*}%
where $\left\langle .,.\right\rangle $ is a pairing of $Lip_{0}\left(
X;Y^{\ast }\right) $ and $\mathcal{F}\left( X;Y\right) $ defined by%
\begin{equation*}
\left\langle T,m\right\rangle =\sum_{x\in X}\left\langle T\left( x\right)
,m\left( x\right) \right\rangle .
\end{equation*}%
Therefore, for a general molecule $m=\sum_{i=1}^{n}\lambda
_{i}y_{i}m_{x_{i}x_{i}^{\prime }},$%
\begin{equation}
\left\langle T,m\right\rangle =\sum_{i=1}^{n}\lambda _{i}\left\langle
T\left( x_{i}\right) -T\left( x_{i}^{\prime }\right) ,y_{i}\right\rangle . 
\tag{0.1}
\end{equation}%
Let $X$ be a pointed metric space and $Y$ be a Banach space, note that for
any $T\in Lip_{0}\left( X;Y\right) ,$ then there exists a unique linear map
(linearization of $T$) $\widehat{T}:\mathcal{F}\left( X\right)
\longrightarrow Y$ such that $\widehat{T}\delta _{X}=T$ and $\left\Vert 
\widehat{T}\right\Vert =Lip\left( T\right) ,$ i.e., the following diagram
commutes%
\begin{equation}
\begin{array}{ccc}
X & \overset{T}{\longrightarrow } & Y \\ 
\delta _{X}\downarrow  & \nearrow \widehat{T} &  \\ 
\mathcal{F}\left( X\right)  &  & 
\end{array}
\tag{0.2}
\end{equation}%
where $\delta _{X}$ is the canonical embedding so that $\left\langle \delta
_{X}\left( x\right) ,f\right\rangle =\left\langle m_{x0},f\right\rangle
=f\left( x\right) $ for $f\in X^{\#}.$ Let $Y$ be a Banach space, then $B_{Y}
$ denotes its closed unit ball and $Y^{\ast }$ its (topological) dual.
Consider $1\leq p\leq \infty $ and $n\in \mathbb{N}^{\ast }$. We denote by $%
l_{p}^{n}\left( Y\right) $ the Banach space of all sequences $\left(
y_{i}\right) _{i=1}^{n}$ in $Y$ with the norm

\begin{center}
$\left\Vert \left( y_{i}\right) _{i}\right\Vert _{l_{p}^{n}\left( Y\right)
}=(\sum_{i=1}^{n}\left\Vert y_{i}\right\Vert ^{p})^{\frac{1}{p}}$,
\end{center}

\noindent and by $l_{p}^{n,\omega }\left( Y\right) $ the Banach space of all
sequences $\left( y_{i}\right) _{i=1}^{n}$ in $Y$ with the norm

\begin{center}
$\left\Vert \left( y_{i}\right) _{i}\right\Vert _{l_{p}^{n,\omega }\left(
Y\right) }=\underset{y^{\ast }\in B_{Y^{\ast }}}{\sup }(\sum_{i=1}^{n}\left%
\vert \left\langle y_{i},y^{\ast }\right\rangle \right\vert ^{p})^{\frac{1}{p%
}}.$
\end{center}

\noindent We also have%
\begin{equation}
\left\Vert \left( y_{i}\right) _{i}\right\Vert _{l_{p}^{n,\omega }\left(
Y\right) }=\underset{y^{\ast }\in B_{Y^{\ast }}}{\sup }(\sum_{i=1}^{n}\left%
\vert \left\langle y_{i},y^{\ast }\right\rangle \right\vert ^{p})^{\frac{1}{p%
}}=\left\Vert \left( y_{i}\right) _{i}\right\Vert _{l_{p}^{n,\omega }\left(
Y^{\ast \ast }\right) }.  \tag{0.3}
\end{equation}

\noindent If $Y=\mathbb{R},$ we simply write $l_{p}^{n}$ and $l_{p}^{n,w}.$

\section{\textsc{The Chevet Saphar norms on the space of molecules}}

Let $E,F$ be Banach spaces, in $\left[ \text{3, 15}\right] $, the
Chevet-Saphar norms $d_{p}$ and $g_{p}$ are defined on tensor product $%
E\otimes F$ for $1\leq p\leq \infty $ as follow 
\begin{equation*}
d_{p}\left( u\right) =\inf \left\{ \left\Vert \left( x_{i}\right)
_{i}\right\Vert _{l_{p^{\ast }}^{n,w}\left( E\right) }\left\Vert \left(
y_{i}\right) _{i}\right\Vert _{l_{p}^{n}\left( F\right) }\right\} ,
\end{equation*}%
taking the infimum over all representations of $u$ of the form $%
u=\sum_{i=1}^{n}x_{i}\otimes y_{i}\in E\otimes F.$ If we interchange the
roles of the weak and strong norms in $d_{p},$ we obtain the transposed norm%
\begin{equation*}
g_{p}\left( u\right) =\inf \left\{ \left\Vert \left( x_{i}\right)
_{i}\right\Vert _{l_{p}^{n}\left( E\right) }\left\Vert \left( y_{i}\right)
_{i}\right\Vert _{l_{p^{\ast }}^{n,w}\left( F\right) }\right\} .
\end{equation*}%
For every $p,$ we have 
\begin{equation*}
g_{p}=d_{p}^{t},
\end{equation*}%
where $d_{p}^{t}$ is defined as follows%
\begin{equation*}
d_{p}^{t}\left( u;E\otimes F\right) =d_{p}\left( u^{t};F\otimes E\right) ,
\end{equation*}%
and the transpose, $u^{t}$, of $u=\sum_{i=1}^{n}x_{i}\otimes y_{i}$ is given
by $u^{t}=\sum_{i=1}^{n}y_{i}\otimes x_{i}.$ Inspired by the tensor norm $%
g_{p},$ we give a new norm on $\mathcal{F}\left( X;Y\right) $ like the one
given by J.A. Chavez Dominguez in $\left[ \text{2}\right] $ for the norm $%
d_{p}$. Note that the space $\mathcal{F}\left( X;Y\right) $ plays the role
of the tensor product in the linear theory. Let $p\in \left[ 1,\infty \right]
$ and $m\in \mathcal{F}\left( X;Y\right) .$ We consider for $m\in \mathcal{F}%
\left( X;Y\right) $%
\begin{equation*}
\mu _{p}\left( m\right) =\inf \left\{ \left\Vert \left( \lambda _{i}d\left(
x_{i},x_{i}^{\prime }\right) \right) _{i}\right\Vert _{l_{p}^{n}}\left\Vert
\left( y_{i}\right) _{i}\right\Vert _{l_{p^{\ast }}^{n,w}\left( Y\right)
}\right\} ,
\end{equation*}%
where the infimum is taken over all representations of $m$ of the form%
\begin{equation*}
m=\sum_{i=1}^{n}\lambda _{i}y_{i}m_{x_{i}x_{i}^{\prime }},
\end{equation*}%
with $x_{i},x_{i}^{\prime }\in X$, $y_{i}\in Y,\lambda _{i}\in \mathbb{R}%
;\left( 1\leq i\leq n\right) $ and $n\in \mathbb{N}^{\ast }.$%
\vspace{0.5cm}%

\textbf{Proposition 1.1.}\textit{\ Let }$X$\textit{\ be a pointed metric
space and }$Y$\textit{\ be a Banach space. Let }$p\in \left[ 1,\infty \right]
,$\textit{\ then }$\mu _{p}$\textit{\ is a norm on }$\mathcal{F}\left(
X;Y\right) .$%
\vspace{0.5cm}%

\textbf{Proof}. It is clear that for any molecule $m\in \mathcal{F}\left(
X;Y\right) $ and any scalar $\alpha $ we have%
\begin{equation*}
\mu _{p}\left( m\right) \geq 0\text{ and }\mu _{p}\left( \alpha m\right)
=\left\vert \alpha \right\vert \mu _{p}\left( m\right) .
\end{equation*}%
Let $y^{\ast }\in Y^{\ast },$ $f\in X^{\#}$ and $m\in \mathcal{F}\left(
X;Y\right) .$ Using the pairing formula (0.1)%
\begin{eqnarray*}
\left\vert \left\langle y^{\ast }f,m\right\rangle \right\vert  &=&\left\vert
\sum_{i=1}^{n}\lambda _{i}y^{\ast }\left( y_{i}\right) \left( f\left(
x_{i}\right) -f\left( x_{i}^{\prime }\right) \right) \right\vert  \\
&\leq &(\sum_{i=1}^{n}\left\vert \lambda _{i}\left( f\left( x_{i}\right)
-f\left( x_{i}^{\prime }\right) \right) \right\vert ^{p})^{\frac{1}{p}%
}(\sum_{i=1}^{n}\left\vert y^{\ast }\left( y_{i}\right) \right\vert
^{p^{\ast }})^{\frac{1}{p^{\ast }}} \\
&\leq &\left\Vert y^{\ast }\right\Vert Lip\left( f\right) \left\Vert \left(
\lambda _{i}d\left( x_{i},x_{i}^{\prime }\right) \right) _{i}\right\Vert
_{l_{p}^{n}}\left\Vert \left( y_{i}\right) _{i}\right\Vert _{l_{p^{\ast
}}^{n,w}\left( Y\right) }.
\end{eqnarray*}%
By taking the infimum over all representations of $m,$ we obtain%
\begin{equation*}
\left\vert \left\langle y^{\ast }f,m\right\rangle \right\vert \leq
\left\Vert y^{\ast }\right\Vert Lip\left( f\right) \mu _{p}\left( m\right) .
\end{equation*}%
Now, suppose that $\mu _{p}\left( m\right) =0$, then for every $y^{\ast }\in
Y^{\ast }$ and $f\in X^{\#}$%
\begin{equation*}
0=\left\langle y^{\ast }f,m\right\rangle =\sum_{i=1}^{n}\left\langle
f,\lambda _{i}y^{\ast }\left( y_{i}\right) m_{x_{i}x_{i}^{\prime
}}\right\rangle ,
\end{equation*}%
by the duality between $\mathcal{F}\left( X\right) $ and $X^{\#}$, the
real-valued molecule $y^{\ast }\circ m$ is equal to $0$ for all $y^{\ast
}\in Y^{\ast }$ and consequently $m=0.$ Let now $m_{1},m_{2}\in \mathcal{F}%
\left( X;Y\right) .$ By the definition of $\mu _{p}$ we can find a
representation%
\begin{equation*}
m_{1}=\dsum\limits_{i=1}^{l}\lambda _{1i}y_{1i}m_{x_{1i}x_{1i}^{\prime }},
\end{equation*}%
such that%
\begin{equation*}
\left\Vert \left( \lambda _{1i}d\left( x_{1i},x_{1i}^{\prime }\right)
\right) _{i}\right\Vert _{l_{p}^{l}}\left\Vert \left( y_{1i}\right)
_{i}\right\Vert _{l_{p^{\ast }}^{l,w}\left( Y\right) }\leq \mu _{p}\left(
m_{1}\right) +\varepsilon .
\end{equation*}%
Replacing $\left( \lambda _{1i}\right) $ and $\left( y_{1i}\right) $ by an
appropriate multiple of them, 
\begin{equation*}
\lambda _{1i}=\lambda _{1i}\frac{\left\Vert \left( y_{1i}\right)
_{i}\right\Vert _{l_{p^{\ast }}^{l,w}\left( Y\right) }^{\frac{1}{p}}}{%
\left\Vert \left( \lambda _{1i}d\left( x_{1i},x_{1i}^{\prime }\right)
\right) _{i}\right\Vert _{l_{p}^{l}}^{\frac{1}{p^{\ast }}}},y_{1i}=y_{1i}%
\frac{\left\Vert \left( \lambda _{1i}d\left( x_{1i},x_{1i}^{\prime }\right)
\right) _{i}\right\Vert _{l_{p}^{l}}^{\frac{1}{p^{\ast }}}}{\left\Vert
\left( y_{1i}\right) _{i}\right\Vert _{l_{p^{\ast }}^{l,w}\left( Y\right) }^{%
\frac{1}{p}}},
\end{equation*}%
we can find%
\begin{equation*}
\left\Vert \left( \lambda _{1i}d\left( x_{1i},x_{1i}^{\prime }\right)
\right) _{i}\right\Vert _{l_{p}^{l}}\leq \left( \mu _{p}\left( m_{1}\right)
+\varepsilon \right) ^{\frac{1}{p}},\text{ }\left\Vert \left( y_{1i}\right)
_{i}\right\Vert _{l_{p^{\ast }}^{l,w}\left( Y\right) }\leq \left( \mu
_{p}\left( m_{1}\right) +\varepsilon \right) ^{\frac{1}{p^{\ast }}}.
\end{equation*}%
Similarly for $m_{2},$ we choose a representation%
\begin{equation*}
m_{2}=\dsum\limits_{i=1}^{s}\lambda _{2i}y_{2i}m_{x_{2i}x_{2i}^{\prime }},
\end{equation*}%
such that%
\begin{equation*}
\left\Vert \lambda _{2i}d\left( x_{2i},x_{2i}^{\prime }\right) \right\Vert
_{l_{p}^{s}}\left\Vert \left( y_{2i}\right) _{i}\right\Vert _{l_{p^{\ast
}}^{s,w}\left( Y\right) }\leq \mu _{p}\left( m_{2}\right) +\varepsilon .
\end{equation*}%
Again, replacing $\left( \lambda _{2i}\right) $ and $\left( y_{2i}\right) $
by an appropriate multiple of them as in the above, we find%
\begin{equation*}
\left\Vert \lambda _{2i}d\left( x_{2i},x_{2i}^{\prime }\right) \right\Vert
_{l_{p}^{s}}\leq \left( \mu _{p}\left( m_{2}\right) +\varepsilon \right) ^{%
\frac{1}{p}},\text{ }\left\Vert \left( y_{2i}\right) _{i}\right\Vert
_{l_{p^{\ast }}^{s,w}\left( Y\right) }\leq \left( \mu _{p}\left(
m_{2}\right) +\varepsilon \right) ^{\frac{1}{p^{\ast }}}.
\end{equation*}%
Now, we have%
\begin{eqnarray*}
&&w_{p}\left( m_{1}+m_{2}\right)  \\
&\leq &(\left\Vert \lambda _{1i}d\left( x_{1i},x_{1i}^{\prime }\right)
\right\Vert _{l_{p}^{l}}^{p}+\left\Vert \lambda _{2i}d\left(
x_{2i},x_{2i}^{\prime }\right) \right\Vert _{l_{p}^{s}}^{p})^{\frac{1}{p}%
}(\left\Vert \left( y_{1i}\right) _{i}\right\Vert _{l_{p^{\ast
}}^{l,w}\left( Y\right) }^{p^{\ast }}+\left\Vert \left( y_{2i}\right)
_{i}\right\Vert _{l_{p^{\ast }}^{s,w}\left( Y\right) }^{p^{\ast }})^{\frac{1%
}{p^{\ast }}} \\
&\leq &\left( \mu _{p}\left( m_{1}\right) +\mu _{p}\left( m_{2}\right)
+2\varepsilon \right) ^{\frac{1}{p}}\left( \mu _{p}\left( m_{1}\right) +\mu
_{p}\left( m_{2}\right) +2\varepsilon \right) ^{\frac{1}{p^{\ast }}} \\
&\leq &\mu _{p}\left( m_{1}\right) +\mu _{p}\left( m_{2}\right)
+2\varepsilon .
\end{eqnarray*}%
By letting $\varepsilon $ tend to zero we obtain the triangle inequality for 
$\mu _{p}.\quad \blacksquare $%
\vspace{0.5cm}%

We denote by $\mathcal{F}_{\mu _{p}}\left( X;Y\right) $ the completion of $%
\mathcal{F}\left( X;Y\right) $ for the norm $\mu _{p}$.%
\vspace{0.5cm}%

\textbf{Proposition 1.2}. \textit{Let }$X$\textit{\ be a pointed metric
space and }$Y$\textit{\ be a Banach space. We have }%
\begin{equation*}
\mathcal{F}_{\mu _{p}}\left( X;Y\right) =\mathcal{F}\left( X\right) \hat{%
\otimes}_{g_{p}}Y,
\end{equation*}%
\textit{where }$g_{p}$\textit{\ is the Chevet-Saphar norm defined as above.}%
\vspace{0.5cm}%

\textbf{Proof}. We can establish the identification via the next linear map

\begin{equation*}
\varphi \left( m\right) =\varphi (\sum_{i=1}^{n}\lambda
_{i}y_{i}m_{x_{i}x_{i}^{\prime }})=u=\sum_{i=1}^{n}\left( \lambda
_{i}m_{x_{i}x_{i}^{\prime }}\right) \otimes y_{i}.
\end{equation*}%
Indeed, we have

\QTP{Body Math}
$\mu _{p}\left( m\right) =\inf \left\{ \left\Vert \left( \lambda _{i}d\left(
x_{i},x_{i}^{\prime }\right) \right) _{i}\right\Vert _{l_{p}^{n}}\left\Vert
\left( y_{i}\right) _{i}\right\Vert _{l_{p^{\ast }}^{n,w}\left( Y\right)
}\right\} $

\QTP{Body Math}
$=\inf \left\{ (\sum_{i=1}^{n}\left\Vert \lambda _{i}m_{x_{i}x_{i}^{\prime
}}\right\Vert _{\mathcal{F}\left( X\right) }^{p})^{\frac{1}{p}}\left\Vert
\left( y_{i}\right) _{i}\right\Vert _{l_{p^{\ast }}^{n,w}\left( Y\right)
}\right\} $

\QTP{Body Math}
$=\inf \left\{ \left\Vert \left( \lambda _{i}m_{x_{i}x_{i}^{\prime }}\right)
_{i}\right\Vert _{l_{p}^{n}\left( \mathcal{F}\left( X\right) \right)
}\left\Vert \left( y_{i}\right) _{i}\right\Vert _{l_{p^{\ast }}^{n,w}\left(
Y\right) }\right\} =g_{p}\left( u\right) .$

\noindent Now, it suffices to show that $\varphi $ is onto. Let $%
u=\sum_{i=1}^{n}v_{i}\otimes y_{i}\in \mathcal{F}\left( X\right) \hat{\otimes%
}_{g_{p}}Y$ where $v_{i}=\dsum\limits_{s_{i}=1}^{k_{i}}\lambda
_{s_{i}}m_{x_{s_{i}}x_{s_{i}}^{\prime }}.$ We put for $1\leq i\leq n,$ $%
m_{i}=\dsum\limits_{s_{i}=1}^{k_{i}}\lambda
_{s_{i}}y_{i}m_{x_{s_{i}}x_{s_{i}}^{\prime }}$ and $m=\sum_{i=1}^{n}m_{i}.$
We will verify that $\varphi \left( m\right) =u$. Indeed,

\QTP{Body Math}
$\varphi (m)=\sum_{i=1}^{n}\varphi \left( m_{i}\right)
=\sum_{i=1}^{n}\varphi (\dsum\limits_{s_{i}=1}^{k_{i}}\lambda
_{s_{i}}y_{i}m_{x_{s_{i}}x_{s_{i}}^{\prime }})$

\QTP{Body Math}
$=\sum_{i=1}^{n}\dsum\limits_{s_{i}=1}^{k_{i}}\lambda
_{s_{i}}m_{x_{s_{i}}x_{s_{i}}^{\prime }}\otimes
y_{i}=\sum_{i=1}^{n}v_{i}\otimes y_{i}=u.\quad \blacksquare 
\vspace{0.5cm}%
$

In the next result we give a characterization of an element of the dual of
the space $\mathcal{F}_{\mu _{p}}\left( X;Y\right) .$\ For the proof, we
need the following lemma which is due to Ky Fan (see $\left[ \text{5},\text{
p 190}\right] $ for more detail about this lemma).%
\vspace{0.5cm}%

\textbf{Lemma} \textbf{1.3} \textbf{(Ky Fan).} \textit{Let }$E$\textit{\ be
a Hausdorff topological vector space and }$\mathcal{C}$\textit{\ be a
compact convex subset of }$E$\textit{. Let }$M$\textit{\ be a set of
functions on} $\mathcal{C}$ \textit{with values in }$(-\infty ,\infty ]$%
\textit{\ having the following:}

\noindent (a)\textit{\ each }$f\in M$\textit{\ is convex and lower
semicontinuous;}

\noindent (b)\textit{\ if }$g\in conv(M)$\textit{, there is an }$f\in M$%
\textit{\ such that }$g(x)\leq f(x)$\textit{, for every} $x\in \mathcal{C}$%
\textit{;}

\noindent (c)\textit{\ there is an }$r\in \mathbb{R}$\textit{\ such that
each }$f\in M$\textit{\ has a value not greater than }$r$\textit{.}

\noindent \textit{Then, there is an} $x_{0}\in \mathcal{C}$ \textit{such
that }$f(x_{0})\leq r$\textit{\ for all }$f\in M$\textit{.}%
\vspace{0.5cm}%

\textbf{Theorem 1.4.} \textit{Let }$X$\textit{\ be a pointed metric space, }$%
Y$\textit{\ be a Banach space and }$C>0.$\textit{\ The following properties
are equivalent.}

\noindent $\left( 1\right) $ \textit{The function }$\varphi $\textit{\ is }$%
\mu _{p}$\textit{-continuous on} $\mathcal{F}\left( X;Y\right) $, \textit{%
i.e.,}%
\begin{equation}
\left\vert \varphi \left( m\right) \right\vert \leq C\mu _{p}\left( m\right) 
\text{ \textit{for all} }m\in \mathcal{F}\left( X;Y\right) .  \tag{1.1}
\end{equation}

\noindent $\left( 2\right) $ \textit{For any representation of} $m$ \textit{%
of the form} $m=\sum_{i=1}^{n}\lambda _{i}y_{i}m_{x_{i}x_{i}^{\prime }}$, 
\textit{we have}%
\begin{equation}
\sum_{i=1}^{n}\left\vert \varphi \left( \lambda
_{i}y_{i}m_{x_{i}x_{i}^{\prime }}\right) \right\vert \leq C\mu _{p}\left(
m\right) .  \tag{1.2}
\end{equation}

\noindent $\left( 3\right) $ \textit{There exists a Radon} \textit{%
probability }$\mu $\textit{\ on }$B_{Y^{\ast }}$ \textit{such that for every
atom of the form }$ym_{xx^{\prime }}$%
\begin{equation}
\left\vert \varphi \left( ym_{xx^{\prime }}\right) \right\vert \leq Cd\left(
x,x^{\prime }\right) \left\Vert y\right\Vert _{L_{p^{\ast }}\left( \mu
\right) }.%
\vspace{0.5cm}
\tag{1.3}
\end{equation}

\textbf{Proof}\textit{. }$\left( 1\right) \implies \left( 2\right) :$ Let $%
\left( \alpha _{i}\right) _{1\leq i\leq n}$ be a scalar sequence. By $\left(
1.1\right) $, we have%
\begin{eqnarray*}
&&\left\vert \varphi (\sum_{i=1}^{n}\alpha _{i}\lambda
_{i}y_{i}m_{x_{i}x_{i}^{\prime }})\right\vert =\left\vert
\sum_{i=1}^{n}\alpha _{i}\varphi \left( \lambda
_{i}y_{i}m_{x_{i}x_{i}^{\prime }}\right) \right\vert  \\
&\leq &C\left\Vert \left( \alpha _{i}\lambda _{i}d\left( x_{i},x_{i}^{\prime
}\right) \right) _{i}\right\Vert _{l_{p}^{n}}\left\Vert \left( y_{i}\right)
_{i}\right\Vert _{l_{p^{\ast }}^{n,w}\left( Y\right) } \\
&\leq &C\left\Vert \left( \alpha _{i}\right) _{i}\right\Vert _{l_{\infty
}^{n}}\left\Vert \left( \lambda _{i}d\left( x_{i},x_{i}^{\prime }\right)
\right) _{i}\right\Vert _{l_{p}^{n}}\left\Vert \left( y_{i}\right)
_{i}\right\Vert _{l_{p^{\ast }}^{n,w}\left( Y\right) }.
\end{eqnarray*}%
Taking the supremum over all sequences $\left( \alpha _{i}\right) _{1\leq
i\leq n}$ with $\left\Vert \left( \alpha _{i}\right) _{i}\right\Vert
_{l_{\infty }^{n}}\leq 1,$ we obtain $\left( 1.2\right) .$

$\left( 2\right) \implies \left( 3\right) :$ Let $\varphi $ be a $\mu _{p}$%
-continuous function on $\mathcal{F}\left( X;Y\right) $. Let $K=B_{Y^{\ast
}} $. We consider the set $\mathcal{C}$ of probability measures on $K$. It
is a convex and compact subset of $C(K)^{\ast }$ endowed with its weak $\ast 
$-topology. Let $M$\ be the set of all functions on $\mathcal{C}$ with
values in $\mathbb{R}$ of the form

\begin{center}
$%
\begin{array}{ll}
& \Psi _{\left( \left( \lambda _{i}\right) ,\left( x_{i}\right) ,\left(
x_{i}^{\prime }\right) ,\left( y_{i}\right) \right) }(\mu ) \\ 
= & \sum_{i=1}^{n}\left\vert \varphi \left( \lambda
_{i}y_{i}m_{x_{i}x_{i}^{\prime }}\right) \right\vert -\sum_{i=1}^{n}(\tfrac{C%
}{p}\left\Vert \left( \lambda _{i}d\left( x_{i},x_{i}^{\prime }\right)
\right) _{i}\right\Vert _{l_{p}^{n}}^{p}+\tfrac{C}{p^{\ast }}\left\Vert
y_{i}\right\Vert _{L_{p^{\ast }}\left( \mu \right) }^{p^{\ast }}),%
\end{array}%
$
\end{center}

\noindent where $\left( x_{i}\right) _{1\leq i\leq n},\left( x_{i}^{\prime
}\right) _{1\leq i\leq n}\subset X,$ $\left( y_{i}\right) _{1\leq i\leq
n}\subset Y$ and $\left( \lambda _{i}\right) _{1\leq i\leq n}\subset \mathbb{%
R}$. We will verify the assumptions of Ky Fan's lemma:\newline
(a)\ It is easy to see that each element of $M$ is convex and continuous on $%
\mathcal{C}.$\newline
(b) It suffices to show that $M$ is convex. Let $\Psi _{1},\Psi _{2}$ in $M$
such that

\begin{center}
$%
\begin{array}{ll}
& \Psi _{_{1}\left( \left( \lambda _{1i}\right) ,\left( x_{1i}\right)
,\left( x_{1i}^{\prime }\right) ,\left( y_{1i}\right) \right) }(\mu ) \\ 
= & \sum\limits_{i=1}^{l}\left\vert \varphi \left( \lambda
_{1i}y_{1i}m_{x_{1i}x_{1i}^{\prime }}\right) \right\vert
-\dsum\limits_{i=1}^{l}(\tfrac{C}{p}\left\Vert \left( \lambda _{1i}d\left(
x_{1i},x_{1i}^{\prime }\right) \right) _{i}\right\Vert _{l_{p}^{l}}^{p}+%
\tfrac{C}{p^{\ast }}\left\Vert y_{1i}\right\Vert _{L_{p^{\ast }}\left( \mu
\right) }^{p^{\ast }}),%
\end{array}%
$
\end{center}

\noindent and

\begin{center}
$%
\begin{array}{ll}
& \Psi _{_{2}\left( \left( \lambda _{2i}\right) ,\left( x_{2i}\right)
,\left( x_{2i}^{\prime }\right) ,\left( y_{2i}\right) \right) }(\mu ) \\ 
= & \sum\limits_{i=1}^{s}\left\vert \varphi \left( \lambda
_{2i}y_{2i}m_{x_{2i}x_{2i}^{\prime }}\right) \right\vert
-\dsum\limits_{i=1}^{s}(\tfrac{C}{p}\left\Vert \left( \lambda _{2i}d\left(
x_{2i},x_{2i}^{\prime }\right) \right) _{i}\right\Vert _{l_{p}^{s}}^{p}+%
\tfrac{C}{p^{\ast }}\left\Vert y_{2i}\right\Vert _{L_{p^{\ast }}\left( \mu
\right) }^{p^{\ast }}).%
\end{array}%
$
\end{center}

\noindent It follows that

\begin{center}
$%
\begin{array}{ll}
& \alpha \Psi _{1}+\left( 1-\alpha \right) \Psi _{2} \\ 
= & \sum_{i=1}^{n}\left\vert \varphi \left( \lambda
_{i}y_{i}m_{x_{i}x_{i}^{\prime }}\right) \right\vert -\sum_{i=1}^{n}(\tfrac{C%
}{p}\left\Vert \left( \lambda _{i}d\left( x_{i},x_{i}^{\prime }\right)
\right) _{i}\right\Vert _{l_{p}^{n}}^{p}+\tfrac{C}{p^{\ast }}\left\Vert
y_{i}\right\Vert _{L_{p^{\ast }}\left( \mu \right) }^{p^{\ast }}),%
\end{array}%
$
\end{center}

\noindent with $n=l+s$, and%
\begin{eqnarray*}
x_{i} &=&\left\{ 
\begin{tabular}{lll}
$x_{1i}$ & if & $1\leq i\leq l$ \\ 
$x_{2\left( i-l\right) }$ & if & $l+1\leq i\leq n$%
\end{tabular}%
\right. , \\
x_{i}^{\prime } &=&\left\{ 
\begin{tabular}{lll}
$x_{1i}^{\prime }$ & if & $1\leq i\leq l$ \\ 
$x_{2\left( i-l\right) }^{\prime }$ & if & $l+1\leq i\leq n$%
\end{tabular}%
\right. , \\
y_{i} &=&\left\{ 
\begin{tabular}{lll}
$\alpha ^{\frac{1}{p^{\ast }}}y_{1i}$ & if & $1\leq i\leq l$ \\ 
$\left( 1-\alpha \right) ^{\frac{1}{p^{\ast }}}y_{2\left( i-l\right) }$ & if
& $l+1\leq i\leq n$%
\end{tabular}%
\right. , \\
\lambda _{i} &=&\left\{ 
\begin{tabular}{lll}
$\alpha ^{\frac{1}{p}}\lambda _{1i}$ & if & $1\leq i\leq l$ \\ 
$\left( 1-\alpha \right) ^{\frac{1}{p}}\lambda _{2\left( i-l\right) }$ & if
& $l+1\leq i\leq n$%
\end{tabular}%
\right. .
\end{eqnarray*}%
(c) Let us show that $r=0$ verifies the condition (c) of Ky Fan's Lemma.
There exists $y_{0}^{\ast }\in B_{Y^{\ast }}$ such that

\begin{center}
$\underset{\left\Vert y^{\ast }\right\Vert _{Y^{\ast }}=1}{\sup }%
(\sum_{i=1}^{n}\left\vert \left\langle y_{i},y^{\ast }\right\rangle
\right\vert ^{p^{\ast }})^{\frac{1}{p^{\ast }}}=(\sum_{i=1}^{n}\left\vert
\left\langle y_{i},y_{0}^{\ast }\right\rangle \right\vert ^{p^{\ast }})^{%
\frac{1}{p^{\ast }}}.$
\end{center}

\noindent Let $\delta _{y_{0}^{\ast }}$ be the Dirac's measure supported by $%
y_{0}^{\ast }$. Using the elementary identity 
\begin{equation*}
\forall \alpha ,\beta \in \mathbb{R}_{+}^{\ast }:\alpha \beta
=\inf_{\epsilon >0}\left\{ \frac{1}{p}\left( \frac{\alpha }{\epsilon }%
\right) ^{p}+\frac{1}{p^{\ast }}\left( \epsilon \beta \right) ^{p^{\ast
}}\right\} ,
\end{equation*}%
and by (1.2), we find that by taking 
\begin{equation*}
\alpha =\left\Vert \left( \lambda _{i}d\left( x_{i},x_{i}^{\prime }\right)
\right) _{i}\right\Vert _{l_{p}^{n}},\beta =(\sum_{i=1}^{n}\left\vert
\left\langle y_{i},y_{0}^{\ast }\right\rangle \right\vert ^{p^{\ast }})^{%
\frac{1}{p}},\epsilon =1.
\end{equation*}

\begin{center}
$%
\begin{array}{ll}
& \Psi \left( \delta _{y_{0}^{\ast }}\right)  \\ 
= & \sum_{i=1}^{n}\left\vert \varphi \left( \lambda
_{i}y_{i}m_{x_{i}x_{i}^{\prime }}\right) \right\vert -\sum_{i=1}^{n}(\tfrac{C%
}{p}\left\Vert \left( \lambda _{i}d\left( x_{i},x_{i}^{\prime }\right)
\right) _{i}\right\Vert _{l_{p}^{n}}^{p}+\tfrac{C}{p^{\ast }}\left\Vert
y_{i}\right\Vert _{L_{p^{\ast }}\left( \delta _{y_{0}^{\ast }}\right)
}^{p^{\ast }}) \\ 
= & \sum_{i=1}^{n}\left\vert \varphi \left( \lambda
_{i}y_{i}m_{x_{i}x_{i}^{\prime }}\right) \right\vert -\sum_{i=1}^{n}(\tfrac{C%
}{p}\left\Vert \left( \lambda _{i}d\left( x_{i},x_{i}^{\prime }\right)
\right) _{i}\right\Vert _{l_{p}^{n}}^{p}+\tfrac{C}{p^{\ast }}\left\vert
\left\langle y_{i},y_{0}^{\ast }\right\rangle \right\vert ^{p^{\ast }}) \\ 
\leq  & \sum_{i=1}^{n}\left\vert \varphi \left( \lambda
_{i}y_{i}m_{x_{i}x_{i}^{\prime }}\right) \right\vert -C\left\Vert \left(
\lambda _{i}d\left( x_{i},x_{i}^{\prime }\right) \right) _{i}\right\Vert
_{l_{p}^{n}}(\sum_{i=1}^{n}\left\vert \left\langle y_{i},y_{0}^{\ast
}\right\rangle \right\vert ^{p^{\ast }})^{\frac{1}{p}} \\ 
\leq  & 0.%
\end{array}%
$
\end{center}

\noindent By Ky Fan's Lemma, there is $\mu \in \mathcal{C}$ such that $\Psi
(\mu )\leq 0$ for every $\Psi \in M$. If we consider $\lambda \in \mathbb{R}%
_{+}^{\ast },x,x^{^{\prime }}\in X$ and $y\in Y,$ we obtain

\begin{center}
$%
\begin{array}{ll}
& \Psi (\mu )=\Psi _{\left( \lambda ,x,x^{\prime },y\right) }(\mu ) \\ 
= & \left\vert \varphi \left( \lambda ym_{xx^{\prime }}\right) \right\vert -%
\tfrac{C}{p}\left\vert \lambda \right\vert ^{p}d\left( x,x^{\prime }\right)
^{p}-\tfrac{C}{p^{\ast }}\left\Vert y\right\Vert _{L_{p^{\ast }}\left( \mu
\right) }^{p^{\ast }}\leq 0.%
\end{array}%
$
\end{center}

\noindent Thus%
\begin{equation*}
\left\vert \lambda \right\vert \left\vert \varphi \left( ym_{xx^{\prime
}}\right) \right\vert \leq \tfrac{C}{p}\left\vert \lambda \right\vert
^{p}d\left( x,x^{\prime }\right) ^{p}+\tfrac{C}{p^{\ast }}\left\Vert
y\right\Vert _{L_{p^{\ast }}\left( \mu \right) }^{p^{\ast }}.
\end{equation*}%
Fix $\epsilon >0$. Replacing $\lambda $ by $\frac{1}{\epsilon ^{p\ast }}$ 
\begin{equation*}
\frac{1}{\epsilon ^{p\ast }}\left\vert \varphi \left( ym_{xx^{\prime
}}\right) \right\vert \leq C(\frac{1}{p\epsilon ^{pp\ast }}d\left(
x,x^{\prime }\right) ^{p}+\tfrac{1}{p^{\ast }}\left\Vert y\right\Vert
_{L_{p^{\ast }}\left( \mu \right) }^{p^{\ast }}).
\end{equation*}%
Then%
\begin{eqnarray*}
\left\vert \varphi \left( ym_{xx^{\prime }}\right) \right\vert &\leq &C(%
\frac{1}{p\epsilon ^{p}}d\left( x,x^{\prime }\right) ^{p}+\tfrac{\epsilon
^{p\ast }}{p^{\ast }}\left\Vert y\right\Vert _{L_{p^{\ast }}\left( K,\mu
\right) }^{p\ast }) \\
&\leq &C(\frac{1}{p}(\frac{d\left( x,x^{\prime }\right) }{\epsilon })^{p}+%
\tfrac{1}{p^{\ast }}(\epsilon \left\Vert y\right\Vert _{L_{p^{\ast }}\left(
\mu \right) })^{p\ast }).
\end{eqnarray*}%
We take the infimum over all $\epsilon >0$, we find%
\begin{equation*}
\left\vert \varphi \left( ym_{xx^{\prime }}\right) \right\vert \leq Cd\left(
x,x^{\prime }\right) \left\Vert y\right\Vert _{L_{p^{\ast }}\left( \mu
\right) }.
\end{equation*}

$\left( 3\right) \implies \left( 1\right) :$ Let $m\in \mathcal{F}\left(
X;Y\right) $ such that 
\begin{equation*}
m=\sum_{i=1}^{n}\lambda _{i}y_{i}m_{x_{i}x_{i}^{\prime }}.
\end{equation*}%
By (1.3)

\QTP{Body Math}
$\left\vert \varphi \left( m\right) \right\vert \leq
\sum_{i=1}^{n}\left\vert \varphi \left( \lambda
_{i}y_{i}m_{x_{i}x_{i}^{\prime }}\right) \right\vert $

\QTP{Body Math}
$\leq C\sum_{i=1}^{n}\lambda _{i}d\left( x_{i},x_{i}^{\prime }\right)
\left\Vert y_{i}\right\Vert _{L_{p^{\ast }}\left( \mu \right) }$ $\left( 
\text{by H\"{o}lder inequality}\right) $

\QTP{Body Math}
$\leq C\left\Vert \left( \lambda _{i}d\left( x_{i},x_{i}^{\prime }\right)
\right) _{i}\right\Vert _{l_{p}^{n}}(\sum_{i=1}^{n}\underset{B_{Y^{\ast }}}{%
\int }\left\vert y^{\ast }\left( y_{i}\right) \right\vert ^{p^{\ast }}d\mu
)^{\tfrac{1}{p^{\ast }}}$

\QTP{Body Math}
$\leq C\left\Vert \left( \lambda _{i}d\left( x_{i},x_{i}^{\prime }\right)
\right) _{i}\right\Vert _{l_{p}^{n}}(\underset{B_{Y^{\ast }}}{\int }%
\sum_{i=1}^{n}\left\vert y^{\ast }\left( y_{i}\right) \right\vert ^{p^{\ast
}}d\mu )^{\tfrac{1}{p^{\ast }}}$

\QTP{Body Math}
$\leq C\left\Vert \left( \lambda _{i}d\left( x_{i},x_{i}^{\prime }\right)
\right) _{i}\right\Vert _{l_{p}^{n}}\left\Vert \left( y_{i}\right)
_{i}\right\Vert _{l_{p^{\ast }}^{n,w}\left( Y\right) },$

\noindent as $m$ is arbitrary, we find%
\begin{equation*}
\left\vert \varphi \left( m\right) \right\vert \leq C\mu _{p}\left( m\right)
,
\end{equation*}%
thus$,$ $\varphi $ is $\mu _{p}$-continuous function on $\mathcal{F}\left(
X;Y\right) $.$\quad \blacksquare $

\section{\textsc{The space of Lipschitz-Cohen strongly }$p$\textsc{-summing}}

Pietsch in $\left[ \text{13}\right] $ has published an interesting monograph
on operator ideals. Recently, several generalizations of certain operator
ideals to the space of Lipschitz maps have been investigate. In this
section, we consider the ideal of (Cohen) strongly $p$-summing operators and
we give its extension to Lipschitz mappings. In $\left[ \text{4}\right] $
Cohen has introduced the following concept: A linear operator $%
u:E\rightarrow F$ between Banach spaces is strongly $p$-summing (or Cohen
strongly $p$-summing) ($1<p\leq \infty $) if there is a positive constant $C$
such that for all $n\in \mathbb{N}^{\ast }$, $x_{1},...,x_{n}\in E$ and $%
y_{1}^{\ast },...,y_{n}^{\ast }\in F^{\ast }$, we have%
\begin{equation}
\sum_{i=1}^{n}\left\vert \left\langle u\left( x_{i}\right) ,y_{i}^{\ast
}\right\rangle \right\vert \leq C(\sum_{i=1}^{n}\left\Vert x_{i}\right\Vert
^{p})^{\frac{1}{p}}\left\Vert \left( y_{i}^{\ast }\right) _{i}\right\Vert
_{l_{p^{\ast }}^{n,w}\left( F^{\ast }\right) }.  \tag{2.1}
\end{equation}%
The smallest constant $C,$ which is noted by $d_{p}(u)$, such that the
inequality (2.1) holds, is called the strongly $p$-summing norm on the space 
$\mathcal{D}_{p}(E;F)$ of all Cohen strongly $p$-summing linear operators
from $E$ into $F$, which is a Banach space. If $p=1$, we have $\mathcal{D}%
_{1}(E;F)=\mathcal{B}(E;F),$ the space of all bounded linear operators from $%
E$ to $F.$%
\vspace{0.5cm}%

We give the same definition to the category of Lipschitz mappings.%
\vspace{0.5cm}%

\textbf{Definition 2.1}. Let $X$ be a pointed metric space and $Y$ be a
Banach space. Let $T:X\rightarrow Y$ be a Lipschitz map. $T$ is
Lipschitz-Cohen strongly $p$-summing ($1<p\leq \infty $)\ if there is a
constant $C>0$\ such that for any $n\in \mathbb{N}^{\ast }$, $\left(
x_{i}\right) _{i},\left( x_{i}^{\prime }\right) _{i}$ in $X;\left(
y_{i}^{\ast }\right) _{i}$ in $Y^{\ast }$ and $\left( \lambda _{i}\right)
_{i}$ in $\mathbb{R}_{+}^{\ast }$ $\left( 1\leq i\leq n\right) $, we have%
\begin{equation}
\sum_{i=1}^{n}\lambda _{i}\left\vert \left\langle T\left( x_{i}\right)
-T\left( x_{i}^{\prime }\right) ,y_{i}^{\ast }\right\rangle \right\vert \leq
C\left\Vert \left( \lambda _{i}d\left( x_{i},x_{i}^{\prime }\right) \right)
_{i}\right\Vert _{l_{p}^{n}}\left\Vert \left( y_{i}^{\ast }\right)
_{i}\right\Vert _{l_{p^{\ast }}^{n,w}\left( Y^{\ast }\right) }.  \tag{2.2}
\end{equation}%
We denote by $\mathcal{D}_{p}^{L}\left( X;Y\right) $ the Banach space of all
Lipschitz-Cohen strongly $p$-summing and $d_{p}^{L}\left( T\right) $ its norm%
\begin{equation*}
d_{p}^{L}(T)=\inf \left\{ C>0\text{, }C\text{ verifying (2.2)}\right\} .
\end{equation*}%
As in the linear case, if $p=1$ we have $\mathcal{D}_{1}^{L}\left(
X;Y\right) =Lip_{0}\left( X;Y\right) .$%
\vspace{0.5cm}%

It is easy to show the next Proposition.%
\vspace{0.5cm}%

\textbf{Proposition 2.2}. \textit{Let }$u$\textit{\ be a bounded linear
operator from }$E$\textit{\ into }$F$\textit{\ and }$1\leq p\leq \infty .$%
\textit{\ Then}%
\begin{equation*}
d_{p}\left( u\right) =d_{p}^{L}\left( u\right) .%
\vspace{0.5cm}%
\end{equation*}

In the next result, we give the Pietsch domination theorem for the class of
Lipschitz-Cohen strongly $p$-summing. The linear version of this theorem has
been given by Pietsch $\left[ \text{14}\right] $ for $p$-summing linear
operators.%
\vspace{0.5cm}%

\textbf{Theorem 2.3 (Pietsch's domination).} \textit{Let }$X$\textit{\ be a
pointed metric space and }$Y$\textit{\ be a Banach space. The following
properties are equivalent.}

\noindent $\left( 1\right) $ \textit{The mapping }$T$ \textit{belongs to }$%
\mathcal{D}_{p}^{L}\left( X;Y\right) .$

\noindent $\left( 2\right) $ \textit{For }$\left( x_{i}\right) _{i},\left(
x_{i}^{\prime }\right) _{i}$\textit{\ in }$X;\left( y_{i}^{\ast }\right) _{i}
$\textit{\ in }$Y^{\ast }$\textit{\ and }$\left( \lambda _{i}\right) _{i}$%
\textit{\ in }$\mathbb{R}_{+}^{\ast }$ $\left( 1\leq i\leq n\right) $\textit{%
, we have}%
\begin{equation}
\left\vert \sum_{i=1}^{n}\lambda _{i}\left\langle T\left( x_{i}\right)
-T\left( x_{i}^{\prime }\right) ,y_{i}^{\ast }\right\rangle \right\vert \leq
d_{p}^{L}\left( T\right) \left\Vert \left( \lambda _{i}d\left(
x_{i},x_{i}^{\prime }\right) \right) _{i}\right\Vert _{l_{p}^{n}}\left\Vert
\left( y_{i}^{\ast }\right) _{i}\right\Vert _{l_{p^{\ast }}^{n,w}\left(
Y^{\ast }\right) }.  \tag{2.3}
\end{equation}

\noindent $\left( 3\right) $ \textit{There exist a constant }$C>0$ \textit{%
and} \textit{a} \textit{Radon} \textit{probability }$\mu $\textit{\ on }$%
B_{Y^{\ast \ast }}$ \textit{such that for all }$x,x^{\prime }\in X$ \textit{%
and} $y^{\ast }\in Y^{\ast },$ \textit{we have}%
\begin{equation}
\left\vert \left\langle T\left( x\right) -T\left( x^{\prime }\right)
,y^{\ast }\right\rangle \right\vert \leq Cd\left( x,x^{\prime }\right)
\left\Vert y^{\ast }\right\Vert _{L_{p^{\ast }}\left( \mu \right) }. 
\tag{2.4}
\end{equation}%
\textit{In this case,}%
\begin{equation*}
d_{p}^{L}(T)=\inf \left\{ C>0\text{, }C\text{ \textit{verifying (2.4)}}%
\right\} .%
\vspace{0.5cm}%
\end{equation*}

\textbf{Proof}\textit{.}

\noindent $\left( 1\right) \implies \left( 2\right) :$ Immediate.

\noindent $\left( 2\right) \implies \left( 3\right) :$ Let $T\in \mathcal{D}%
_{p}^{L}\left( X;Y\right) .$ First, we can see $T$ as a mapping defined from 
$X$ into $Y^{\ast \ast }$. Let $\varphi _{T}$ its correspondent linear
function on $\mathcal{F}\left( X;Y^{\ast }\right) $. Let $m\in \mathcal{F}%
\left( X;Y^{\ast }\right) $, by $\left( 2.3\right) $ we have%
\begin{eqnarray*}
\left\vert \varphi _{T}\left( m\right) \right\vert  &=&\left\vert
\sum_{i=1}^{n}\lambda _{i}\left\langle T\left( x_{i}\right) -T\left(
x_{i}^{\prime }\right) ,y_{i}^{\ast }\right\rangle \right\vert  \\
&\leq &C\left\Vert \left( \lambda _{i}d\left( x_{i},x_{i}^{\prime }\right)
\right) _{i}\right\Vert _{l_{p}^{n}}\left\Vert \left( y_{i}^{\ast }\right)
_{i}\right\Vert _{l_{p^{\ast }}^{n,w}\left( Y^{\ast }\right) },
\end{eqnarray*}%
hence, as $m$ is arbitrary, 
\begin{equation*}
\left\vert \varphi _{T}\left( m\right) \right\vert \leq C\mu _{p}\left(
m\right) ,
\end{equation*}%
then $\varphi _{T}$ is $\mu _{p}$-continuous on $\mathcal{F}\left( X;Y^{\ast
}\right) .$ By Theorem 1.4 (3), we have for all $x,x^{\prime }\in X$ and $%
y^{\ast }\in Y^{\ast }$%
\begin{eqnarray*}
\left\vert \left\langle T\left( x\right) -T\left( x^{\prime }\right)
,y^{\ast }\right\rangle \right\vert  &=&\left\vert \varphi _{T}\left(
y^{\ast }m_{xx^{\prime }}\right) \right\vert  \\
&\leq &Cd\left( x,x^{\prime }\right) \left\Vert y^{\ast }\right\Vert
_{L_{p^{\ast }}\left( \mu \right) }.
\end{eqnarray*}

\noindent $\left( 3\right) \implies \left( 1\right) :$ Let $T$ be a
Lipschitz mapping verifies (2.4). For $x,x^{\prime }\in X$ and $y^{\ast }\in
Y^{\ast },$%
\begin{eqnarray*}
\left\vert \varphi _{T}\left( y^{\ast }m_{xx^{\prime }}\right) \right\vert 
&=&\left\vert \left\langle T\left( x\right) -T\left( x^{\prime }\right)
,y^{\ast }\right\rangle \right\vert  \\
&\leq &Cd\left( x,x^{\prime }\right) \left\Vert y^{\ast }\right\Vert
_{L_{p^{\ast }}\left( \mu \right) }.
\end{eqnarray*}%
So, by Theorem 1.4 (1) $\varphi _{T}$ is $\mu _{p}$-continuous on $\mathcal{F%
}\left( X;Y^{\ast }\right) $ and by $\left( 1.2\right) $%
\begin{eqnarray*}
\sum_{i=1}^{n}\lambda _{i}\left\vert \left\langle T\left( x_{i}\right)
-T\left( x_{i}^{\prime }\right) ,y_{i}^{\ast }\right\rangle \right\vert 
&=&\sum_{i=1}^{n}\left\vert \varphi _{T}\left( \lambda _{i}y_{i}^{\ast
}m_{x_{i}x_{i}^{\prime }}\right) \right\vert  \\
&\leq &C\mu _{p}\left( m\right)  \\
&\leq &C\left\Vert \left( \lambda _{i}d\left( x_{i},x_{i}^{\prime }\right)
\right) _{i}\right\Vert _{l_{p}^{n}}\left\Vert \left( y_{i}^{\ast }\right)
_{i}\right\Vert _{l_{p^{\ast }}^{n,w}\left( Y^{\ast }\right) }.
\end{eqnarray*}%
Therefore $T$ is in $\mathcal{D}_{p}^{L}\left( X;Y\right) $ and 
\begin{equation*}
d_{p}^{L}\left( T\right) \leq C.\quad \blacksquare 
\vspace{0.5cm}%
\end{equation*}

The main result of this section is the following identification.%
\vspace{0.5cm}%

\textbf{Theorem 2.4.} \textit{Let }$X$\textit{\ be a pointed metric space
and }$Y$\textit{\ be a Banach space. Let }$p\in \left[ 1,\infty \right] .$ 
\textit{We have the isometric identification}%
\begin{equation}
\mathcal{D}_{p}^{L}\left( X;Y^{\ast }\right) =\mathcal{F}_{\mu
_{_{p}}}\left( X;Y\right) ^{\ast }.%
\vspace{0.5cm}
\tag{2.5}
\end{equation}

\textbf{Proof}. Let $T\in \mathcal{D}_{p}^{L}\left( X;Y^{\ast }\right) $ and 
$\varphi _{T}$ its correspondent linear function on $\mathcal{F}\left(
X;Y\right) .$ We will show that $\varphi _{T}$ is $\mu _{p}$-continuous. Let 
$m=\sum_{i=1}^{n}\lambda _{i}y_{i}m_{x_{i}x_{i}^{\prime }}\in \mathcal{F}%
\left( X;Y\right) $. As $y_{i}$ is an element in $Y^{\ast \ast }$, we obtain
Theorem 2.3

\QTP{Body Math}
$\left\vert \varphi _{T}\left( m\right) \right\vert =\left\vert
\sum_{i=1}^{n}\lambda _{i}\left\langle T\left( x_{i}\right) -T\left(
x_{i}^{\prime }\right) ,y_{i}\right\rangle \right\vert $

\QTP{Body Math}
$\leq d_{p}^{L}\left( T\right) \left\Vert \left( \lambda _{i}d\left(
x_{i},x_{i}^{\prime }\right) \right) _{i}\right\Vert
_{l_{p}^{n}}\sup_{y^{\ast \ast \ast }\in B_{Y^{\ast \ast \ast
}}}(\sum_{i=1}^{n}\left\vert y^{\ast \ast \ast }\left( y_{i}\right)
\right\vert ^{p^{\ast }})^{\frac{1}{p^{\ast }}}$

\QTP{Body Math}
$\leq d_{p}^{L}\left( T\right) \left\Vert \left( \lambda _{i}d\left(
x_{i},x_{i}^{\prime }\right) \right) _{i}\right\Vert _{l_{p}^{m}}\left\Vert
\left( y_{i}\right) _{i}\right\Vert _{l_{p^{\ast }}^{n,w}\left( Y\right) }.$

\noindent Hence, as $m$ is arbitrary, 
\begin{equation*}
\left\vert \varphi _{T}\left( m\right) \right\vert \leq d_{p}^{L}\left(
T\right) \mu _{p}\left( m\right) ,
\end{equation*}%
then $\varphi _{T}$ is $\mu _{p}$-continuous on $\mathcal{F}\left(
X;Y\right) $ and $\left\Vert \varphi _{T}\right\Vert _{\mu _{p}}\leq
d_{p}^{L}\left( T\right) .$

\noindent Conversely, let $\varphi \in \mathcal{F}_{\mu _{p}}\left(
X;Y\right) ^{\ast }$. Note that $\varphi $ can be identified with a mapping $%
T_{\varphi }:X\longrightarrow Y^{\ast }$ via the formula 
\begin{equation*}
\left\langle T_{\varphi }\left( x\right) ,y\right\rangle =\varphi \left(
ym_{x0}\right) .
\end{equation*}%
It is clear that $T_{\varphi }$ is Lipschitz. Indeed, 
\begin{eqnarray*}
\left\Vert T_{\varphi }\left( x\right) -T_{\varphi }\left( x^{\prime
}\right) \right\Vert &=&\sup_{y\in B_{Y}}\left\vert \left\langle T_{\varphi
}\left( x\right) -T_{\varphi }\left( x^{\prime }\right) ,y\right\rangle
\right\vert \\
&=&\sup_{y\in B_{Y}}\left\vert \left\langle \varphi \left( ym_{x0}\right)
-\varphi \left( ym_{x^{\prime }0}\right) \right\rangle \right\vert \\
&=&\sup_{y\in B_{Y}}\left\vert \varphi \left( ym_{xx^{\prime }}\right)
\right\vert \\
&\leq &\sup_{y\in B_{Y}}\left\Vert \varphi \right\Vert d\left( x,x^{\prime
}\right) \left\Vert y\right\Vert \\
&\leq &\left\Vert \varphi \right\Vert d\left( x,x^{\prime }\right) .
\end{eqnarray*}%
Now, let $\left( x_{i}\right) _{i},\left( x_{i}^{\prime }\right) _{i}\subset
X,$ $\left( y_{i}^{\ast \ast }\right) _{i}\subset Y^{\ast \ast }$ and $%
\left( \lambda _{i}\right) _{i}\subset \mathbb{R}_{+}^{\ast }$ $\left( 1\leq
i\leq n\right) .$ Consider the finite-dimensional subspaces 
\begin{equation*}
V=span\left( y_{i}^{\ast \ast }\right) _{i=1}^{n}\subset Y^{\ast \ast },
\end{equation*}%
and $U=span\left( T_{\varphi }\left( x_{i}\right) -T_{\varphi }\left(
x_{i}^{\prime }\right) \right) _{i=1}^{n}\subset Y^{\ast }$. Let $%
\varepsilon >0$. By the principle of local reflexivity $\left[ \text{5, page
178}\right] $, there is an injective linear map $\phi :V\rightarrow Y$ such
that 
\begin{equation*}
\max \left\{ \left\Vert \phi \right\Vert ,\left\Vert \phi \right\Vert
\left\Vert \phi ^{-1}\right\Vert \right\} \leq 1+\varepsilon ,
\end{equation*}%
and $\left\langle \phi \left( y^{\ast \ast }\right) ,u^{\ast }\right\rangle
=\left\langle y^{\ast \ast },u^{\ast }\right\rangle $ for any $y^{\ast \ast
}\in V$ and $u^{\ast }\in U.$ Letting $y_{i}=\phi \left( y_{i}^{\ast \ast
}\right) $, the latter condition together with the continuity of $\varphi $
imply

\QTP{Body Math}
$\left\vert \sum_{i=1}^{n}\lambda _{i}\left\langle T_{\varphi }\left(
x_{i}\right) -T_{\varphi }\left( x_{i}^{\prime }\right) ,y_{i}^{\ast \ast
}\right\rangle \right\vert =\left\vert \sum_{i=1}^{n}\lambda
_{i}\left\langle T_{\varphi }\left( x_{i}\right) -T_{\varphi }\left(
x_{i}^{\prime }\right) ,y_{i}\right\rangle \right\vert $

\QTP{Body Math}
$=\left\vert \varphi (\sum_{i=1}^{n}\lambda _{i}y_{i}m_{x_{i}x_{i}^{\prime
}})\right\vert \leq \left\Vert \varphi \right\Vert _{\mu _{p}}\mu
_{p}(\sum_{i=1}^{n}\lambda _{i}y_{i}m_{x_{i}x_{i}^{\prime }})$

\QTP{Body Math}
$\leq \left\Vert \varphi \right\Vert _{\mu _{p}}\left\Vert \left( \lambda
_{i}d\left( x_{i},x_{i}^{\prime }\right) \right) _{i}\right\Vert
_{l_{p}^{n}}\left\Vert \left( y_{i}\right) _{i}\right\Vert _{l_{p^{\ast
}}^{n,w}\left( Y\right) }.$

\noindent Noting that 
\begin{eqnarray*}
\left\Vert \left( y_{i}\right) _{i}\right\Vert _{l_{p^{\ast }}^{n,w}\left(
Y\right) } &=&\sup_{\left\Vert y^{\ast }\right\Vert \leq
1}(\sum_{i=1}^{n}\left\vert y^{\ast }\left( y_{i}\right) \right\vert
^{p^{\ast }})^{\frac{1}{p^{\ast }}} \\
&=&\sup_{\left\Vert y^{\ast }\right\Vert \leq 1}(\sum_{i=1}^{n}\left\vert
y^{\ast }\left( \phi \left( y_{i}^{\ast \ast }\right) \right) \right\vert
^{p^{\ast }})^{\frac{1}{p^{\ast }}} \\
&\leq &\left\Vert \phi \right\Vert \left\Vert \left( y_{i}^{\ast \ast
}\right) _{i}\right\Vert _{l_{p^{\ast }}^{n,w}\left( Y^{\ast \ast }\right) }
\\
&\leq &\left( 1+\varepsilon \right) \left\Vert \left( y_{i}^{\ast \ast
}\right) _{i}\right\Vert _{l_{p^{\ast }}^{n,w}\left( Y^{\ast \ast }\right) }.
\end{eqnarray*}

\noindent Since $\varepsilon >0$ was arbitrary, letting it go to zero proves
that $T_{\varphi }$ is Lipschitz-Cohen strongly $p$-summing and $%
d_{p}^{L}\left( T_{\varphi }\right) \leq \left\Vert \varphi \right\Vert
_{\mu _{p}}.\quad \blacksquare 
\vspace{0.5cm}%
$

\section{\textsc{Some inclusion and coincidence\ properties}}

The aim of this section is to explore more properties of the class of
Lipschitz-Cohen $p$-summing operators. We start by showing the relationship
between the Lipschitz mapping and its linearization for the concept of
strongly $p$-summing. A similar characterization holds for Lipschitz compact
operators (see $\left[ \text{11}\right] $).%
\vspace{0.5cm}%

\textbf{Proposition 3.1}. \textit{The following properties are equivalent.}

\noindent $\left( 1\right) $ \textit{The mapping }$T$ \textit{belongs to }$%
\mathcal{D}_{p}^{L}\left( X;Y\right) .$

\noindent $\left( 2\right) $ \textit{The linear operator }$\widehat{T}$ 
\textit{belongs to }$\mathcal{D}_{p}\left( \mathcal{F}\left( X\right)
;Y\right) .$

\noindent Even more\textit{, }$\mathcal{D}_{p}^{L}\left( X;Y\right) =%
\mathcal{D}_{p}\left( \mathcal{F}\left( X\right) ;Y\right) $ \textit{holds
isometrically.}%
\vspace{0.5cm}%

\textbf{Proof}. First, suppose that $T\in \mathcal{D}_{p}^{L}\left(
X;Y\right) .$ Let $m\in \mathcal{F}\left( X\right) $ and $y^{\ast }\in
Y^{\ast }$. Then%
\begin{eqnarray*}
\left\vert \left\langle \widehat{T}\left( m\right) ,y^{\ast }\right\rangle
\right\vert  &\leq &\sum_{i=1}^{n}\left\vert \lambda _{i}\right\vert
\left\vert \left\langle T\left( x_{i}\right) -T\left( x_{i}^{\prime }\right)
,y^{\ast }\right\rangle \right\vert  \\
&\leq &d_{p}^{L}\left( T\right) \sum_{i=1}^{n}\left\vert \lambda
_{i}\right\vert d\left( x_{i},x_{i}^{\prime }\right) \left\Vert y^{\ast
}\right\Vert _{L_{p^{\ast }}\left( \mu \right) },
\end{eqnarray*}%
as $m$ is arbitrary, we obtain%
\begin{equation*}
\left\vert \left\langle \widehat{T}\left( m\right) ,y^{\ast }\right\rangle
\right\vert \leq d_{p}^{L}\left( T\right) \left\Vert m\right\Vert _{\mathcal{%
F}\left( X\right) }\left\Vert y^{\ast }\right\Vert _{L_{p^{\ast }}\left( \mu
\right) },
\end{equation*}%
hence $\widehat{T}$ verifies the Pietsch's domination for (Cohen) strongly $p
$-summing operators $\left[ \text{4, Theorem 2.3.1}\right] $, then $\widehat{%
T}\in \mathcal{D}_{p}\left( \mathcal{F}\left( X\right) ;Y\right) $ and 
\begin{equation*}
d_{p}\left( \widehat{T}\right) \leq d_{p}^{L}\left( T\right) .
\end{equation*}%
Conversely, suppose that $\widehat{T}\in \mathcal{D}_{p}\left( \mathcal{F}%
\left( X\right) ;Y\right) .$ Let $x,x^{\prime }\in X$ and $y^{\ast }\in
Y^{\ast },$ by Pietsch's domination of $\widehat{T}$ 
\begin{eqnarray*}
\left\vert \left\langle T\left( x\right) -T\left( x^{\prime }\right)
,y^{\ast }\right\rangle \right\vert  &=&\left\vert \left\langle \widehat{T}%
\left( m_{xx^{\prime }}\right) ,y^{\ast }\right\rangle \right\vert  \\
&\leq &d_{p}\left( \widehat{T}\right) \left\Vert m_{xx^{\prime }}\right\Vert
\left\Vert y^{\ast }\right\Vert _{L_{p^{\ast }}\left( \mu \right) } \\
&\leq &d_{p}\left( \widehat{T}\right) d\left( x,x^{\prime }\right)
\left\Vert y^{\ast }\right\Vert _{L_{p^{\ast }}\left( \mu \right) },
\end{eqnarray*}%
by Theorem 2.3, $T$ is in $\mathcal{D}_{p}^{L}\left( X;Y\right) $ and%
\begin{equation*}
d_{p}^{L}\left( T\right) \leq d_{p}\left( \widehat{T}\right) .\quad
\blacksquare 
\vspace{0.5cm}%
\end{equation*}

One of the nice results of Cohen is that a linear map $u:E\rightarrow F$
between Banach spaces is strongly $p$-summing if and only if the adjoint map 
$u^{\ast }:F^{\ast }\rightarrow E^{\ast }$ is $p^{\ast }$-summing. It would
be interesting to point out that an analogous situation holds in the
nonlinear case: if $X$ is a metric space and $Y$ is a Banach space, $%
T:X\rightarrow Y$ is Lipschitz-Cohen strongly $p$-summing if and only if the
"adjoint" map $T^{\#}\mid _{Y^{\ast }}:Y^{\ast }\rightarrow X^{\#}$ is $%
p^{\ast }$-summing (this map is actually just the linear adjoint of the
linearization $\widehat{T}:\mathcal{F}\left( X\right) \rightarrow Y$).%
\vspace{0.5cm}%

For Lipschitz $p$-summing operators we have the following result.$%
\vspace{0.5cm}%
$

\textbf{Proposition 3.2}. \textit{Let} $1\leq p<\infty $\textit{. Let }$%
T:X\rightarrow Y$ \textit{be a Lipschitz map and }$\widehat{T}$ \textit{its
linearization. Suppose that }$\widehat{T}$\textit{\ is }$p$\textit{-summing,
then }$T$\textit{\ is Lipschitz }$p$\textit{-summing.}%
\vspace{0.5cm}%

\textbf{Proof}. If $\widehat{T}$ is $p$-summing then it is Lipschitz $p$%
-summing, and by (0.2) $T$ Lipschitz factors through $\widehat{T},$ so $T$
is Lipschitz $p$-summing by the ideal property of Lipschitz $p$-summing
operators.$\quad \blacksquare 
\vspace{0.5cm}%
$

\textbf{Remark 3.3}. The converse of the precedent Proposition is not true.
Indeed, the canonical inclusion%
\begin{equation*}
\delta _{\mathbb{R}}:\mathbb{R\rightarrow }\delta _{\mathbb{R}}\left( 
\mathbb{R}\right) \subset \mathcal{F}\left( \mathbb{R}\right)
\end{equation*}%
is Lipschitz $p$-summing since it is Lipschitz equivalent to the identity $%
id_{\mathbb{R}};$ hence $\delta _{\mathbb{R}}:\mathbb{R\rightarrow }\mathcal{%
F}\left( \mathbb{R}\right) $ is Lipschitz $p$-summing. But the linearization
of this map is the identity on $\mathcal{F}\left( \mathbb{R}\right) $, which
cannot be $p$-summing because $\mathcal{F}\left( \mathbb{R}\right) $ is
infinite-dimensional (isometric to $L_{1}\left( \mathbb{R}\right) ,$ in
fact).%
\vspace{0.5cm}%

\textbf{Corollary 3.4}. \textit{Let} $X$\textit{\ be a pointed metric space
and }$Y$ \textit{is an} $\mathcal{L}_{p}$\textit{-space} $\left( 1\leq
p<\infty \right) $\textit{. Then }%
\begin{equation*}
\mathcal{D}_{p^{\ast }}^{L}\left( X;Y\right) \subset \Pi _{p}^{L}\left(
X;Y\right) .
\end{equation*}%
\textbf{Proof}. If $T$ is in $\mathcal{D}_{p^{\ast }}^{L}\left( X;Y\right) $%
, the Proposition 3.1 implies that $\widehat{T}:\mathcal{F}\left( X\right)
\rightarrow Y$ is Cohen strongly $p^{\ast }$-summing. By a result of Cohen $%
\left[ \text{4, Theorem 3.2.3}\right] $, $\widehat{T}$ is $p$-summing and by
Proposition 3.2, $T$ is Lipschitz $p$-summing with 
\begin{equation*}
\pi _{p}^{L}\left( T\right) \leq d_{p^{\ast }}^{L}\left( T\right) .\quad
\blacksquare 
\vspace{0.5cm}%
\end{equation*}

We recall that (see $\left[ \text{2}\right] $) $cs_{p}\left( X;Y\right) $ is
the space of molecules $\mathcal{F}\left( X;Y\right) $ endowed with the next
norm%
\begin{equation*}
cs_{p}\left( m\right) =\inf \left\{ \left\Vert \left( \lambda _{i}\left\Vert
y_{i}\right\Vert \right) _{i}\right\Vert _{l_{p}^{n}}w_{p^{\ast
}}^{Lip}\left( \left( \lambda _{i}^{-1},x_{i},x_{i}^{\prime }\right)
_{i}\right) \right\} ,
\end{equation*}%
where the infimum is taken over all representations of $m$ of the form%
\begin{equation*}
m=\sum_{i=1}^{n}y_{i}m_{x_{i}x_{i}^{\prime }},
\end{equation*}%
with $x_{i},x_{i}^{\prime }\in X$, $y_{i}\in Y,\lambda _{i}\in \mathbb{R}%
_{+}^{\ast };\left( 1\leq i\leq n\right) $ and $n\in \mathbb{N}^{\ast }.$%
\vspace{0.5cm}%

\textbf{Corollary 3.5. }\textit{Let} $X$ \textit{be} \textit{a pointed
metric space, }$1<p<\infty $\textit{\ and }$Y$ \textit{be an} $\mathcal{L}%
_{p}$\textit{-space. The identity mapping}%
\begin{equation*}
id:cs_{p}\left( X;Y\right) \longrightarrow \mathcal{F}_{\mu _{p}}\left(
X;Y\right)
\end{equation*}%
\textit{is continuous with} $\left\Vert id\right\Vert \leq 1.%
\vspace{0.5cm}%
$

\textbf{Proof}. Let $m\in cs_{p}\left( X;Y\right) $ and $\varphi \in 
\mathcal{F}_{\mu _{p}}\left( X;Y\right) ^{\ast }$ such that 
\begin{equation*}
\left\Vert \varphi \right\Vert _{\mu _{p}}\leq 1.
\end{equation*}%
By Theorem 2.4, we can identify $\varphi $ with a function $T_{\varphi }\in 
\mathcal{D}_{p}^{L}\left( X;Y^{\ast }\right) $ with $\left\Vert \varphi
\right\Vert _{\mu _{p}}=d_{p}^{L}\left( T_{\varphi }\right) .$ By the above
result, $T_{\varphi }\in \Pi _{p^{\ast }}^{L}\left( X;Y^{\ast }\right) ,$
and Theorem 4.3 in $\left[ \text{2}\right] $ asserts that $\varphi \in
cs_{p}\left( X;Y\right) ^{\ast }$ with 
\begin{eqnarray*}
\left\Vert \varphi \right\Vert _{cs_{p}\left( X;Y\right) ^{\ast }} &=&\pi
_{p^{\ast }}^{L}\left( T_{\varphi }\right) \\
&\leq &d_{p}^{L}\left( T_{\varphi }\right) =\left\Vert \varphi \right\Vert
_{\mu _{p}},
\end{eqnarray*}%
consequently,%
\begin{eqnarray*}
\mu _{p}\left( m\right) &=&\sup_{\left\Vert \varphi \right\Vert _{\mu
_{p}}\leq 1}\left\Vert \varphi \left( m\right) \right\Vert \\
&\leq &\sup_{\left\Vert \varphi \right\Vert _{cs_{p}\left( X;Y\right) ^{\ast
}}\leq 1}\left\Vert \varphi \left( m\right) \right\Vert =\left\Vert
m\right\Vert _{cs_{p}\left( X;Y\right) }.\quad \blacksquare 
\vspace{0.5cm}%
\end{eqnarray*}

In the next result, we give a version of Grothendieck's Theorem (this famous
result is due to Grothendieck $\left[ \text{10}\right] $). We mention that
other nonlinear versions have already appeared in the literature (for
example in $[$7$]$).%
\vspace{0.5cm}%

\textbf{Corollary 3.6}. (\textit{Grothendieck's\ Theorem}) \textit{Let} $%
X=l_{1}$\textit{\ and }$H$ \textit{be a Hilbert space}.\textit{\ Then }%
\begin{equation*}
\Pi _{1}^{L}\left( X;H\right) =Lip_{0}\left( X;H\right) .%
\vspace{0.5cm}%
\end{equation*}%
\textbf{Proof}. In this case, the free Banach space $\mathcal{F}\left(
X\right) $ is isometrically isomorphic to $L_{1}\left( \mathbb{R}\right) $
(see $\left[ \text{6, Corollary 8}\right] $), then $\widehat{T}:\mathcal{F}%
\left( X\right) \rightarrow H$ is $1$-summing, consequently $T$ is $1$%
-summing.$\quad \blacksquare $%
\vspace{0.5cm}%

In the last result, we consider Lipschitz $\left( p,r,s\right) $-summing
linear operators and we combine with Theorem 5.2 and 5.4 in $\left[ \text{2}%
\right] $ for giving a factorization result using the language of
Lipschitz-Cohen strongly $p$-summing operators.%
\vspace{0.5cm}%

We recall the following definition as stated in $\left[ \text{2}\right] $.

\textbf{Definition 3.7}. Let $X$ be a pointed metric space and $Y$ be a
Banach space. Let $T:X\rightarrow Y$ be a Lipschitz map. $T$ is Lipschitz $%
\left( p,r,s\right) $-summing if there is a constant $C>0$\ such that for
any $n\in \mathbb{N}^{\ast }$, $\left( x_{i}\right) _{i},\left(
x_{i}^{\prime }\right) _{i}$ in $X;\left( y_{i}^{\ast }\right) _{i}$ in $%
Y^{\ast }$ and $\left( \lambda _{i}\right) _{i},\left( k_{i}\right) _{i}$ in 
$\mathbb{R}_{+}^{\ast }$ $\left( 1\leq i\leq n\right) $, we have%
\begin{equation}
\left\Vert \left( \lambda _{i}\left\langle T\left( x_{i}\right) -T\left(
x_{i}^{\prime }\right) ,y_{i}^{\ast }\right\rangle \right) _{i}\right\Vert
_{l_{p}^{n}}\leq Cw_{r}^{Lip}\left( \left( \lambda
_{i}k_{i}^{-1},x_{i},x_{i}^{\prime }\right) _{i}\right) \left\Vert \left(
k_{i}y_{i}^{\ast }\right) _{i}\right\Vert _{l_{s}^{n,w}\left( Y^{\ast
}\right) }.  \tag{3.1}
\end{equation}%
We denote by $\Pi _{p,r,s}^{L}\left( X;Y\right) $ the Banach space of all
Lipschitz $\left( p,r,s\right) $-summing.%
\vspace{0.5cm}%

\textbf{Theorem 3.8}. \textit{Let }$p,r,s\in \left[ 1,\infty \right] $%
\textit{\ such that }$\frac{1}{p}+\frac{1}{r}+\frac{1}{s}=1.$\textit{\ Let }$%
T\in Lip_{0}\left( X;Y\right) $\textit{, the following are equivalent.}

\noindent $\left( 1\right) $\textit{\ The mapping }$T$\textit{\ belongs to }$%
\Pi _{p^{\ast },r,s}^{L}\left( X;Y\right) .$

\noindent $\left( 2\right) $\textit{\ There exist a constant }$C>0$ \textit{%
and} \textit{regular Borel probability measures }$\mu $\textit{\ and }$\nu $%
\textit{\ on the weak}$^{\ast }$\textit{\ compact unit balls }$%
B_{X^{\#}},B_{Y^{\ast \ast }}$\textit{\ such that for all }$x,x^{\prime }\in
X$\textit{\ and }$y^{\ast }\in Y^{\ast }$%
\begin{eqnarray}
&&\left\vert \left\langle T\left( x\right) -T\left( x^{\prime }\right)
,y^{\ast }\right\rangle \right\vert  \TCItag{3.2} \\
&\leq &C(\int_{B_{X^{\#}}}\left\vert f\left( x\right) -f\left( x^{\prime
}\right) \right\vert ^{r}d\mu \left( f\right) )^{\frac{1}{r}%
}(\int_{B_{Y^{\ast \ast }}}\left\vert y^{\ast \ast }\left( y^{\ast }\right)
\right\vert ^{s}d\nu \left( y^{\ast \ast }\right) )^{\frac{1}{s}}.  \notag
\end{eqnarray}

\noindent $\left( 3\right) $\textit{\ There exist a metric space }$\tilde{X}$
\textit{and two Lipschitz mappings }$T_{1},T_{2}$\textit{\ such that }$%
T_{1}\in \Pi _{r}^{L}\left( X;\tilde{X}\right) ,T_{2}\in \mathcal{D}%
_{s^{\ast }}^{L}\left( \tilde{X};Y\right) $\textit{\ and}%
\begin{equation}
T=T_{2}\circ T_{1}.%
\vspace{0.5cm}
\tag{3.3}
\end{equation}

\textbf{Proof}. $\left( 1\right) \iff \left( 2\right) :$ First, by (0.3) the
norm of $T$ as an element of $\Pi _{p^{\ast },r,s}^{L}\left( X;Y\right) $ is
the same as its norm in $\Pi _{p^{\ast },r,s}^{L}\left( X;Y^{\ast \ast
}\right) .$ Then, the equivalence follows from $\left[ \text{2, Theorem 5.2
and 5.4}\right] $ (specialized to the case $E=Y^{\ast }$).

\noindent $\left( 2\right) \implies \left( 3\right) :$ Suppose that $T$
verifies (3.2). Then, we have the following diagram which is commutative

\begin{center}
$%
\begin{tabular}{lll}
$X$ & $\overset{T}{\longrightarrow }$ & $Y$ \\ 
$\downarrow j_{X}$ & $\nearrow \overline{T}$ &  \\ 
$\tilde{X}$ & $\subset L_{r}\left( \mu \right) $ & 
\end{tabular}%
$
\end{center}

\noindent where $j_{X}:X\rightarrow L_{r}\left( \mu \right) $ is the
isometric injection (which is Lipschitz $r$-summing), $\tilde{X}=j_{X}(X)$
is a pointed metric space of which the metric is defined by%
\begin{equation*}
\text{For }\tilde{x},\tilde{x}^{\prime }\in \tilde{X}:d\left( \tilde{x},%
\tilde{x}^{\prime }\right) =\left\Vert \tilde{x}-\tilde{x}^{\prime
}\right\Vert _{L_{r}\left( \mu \right) },
\end{equation*}%
and its origin is $j_{X}\left( 0\right) $. We have $T=\overline{T}\circ
j_{X}.$ The mapping $\overline{T}$ is well defined and is Lipschitz-Cohen
strongly $s^{\ast }$-summing. Indeed,%
\begin{eqnarray*}
\left\vert \left\langle \overline{T}\left( \tilde{x}\right) -\overline{T}%
\left( \tilde{x}^{\prime }\right) ,y^{\ast }\right\rangle \right\vert
&=&\left\vert \left\langle T\left( x\right) -T\left( x^{\prime }\right)
,y^{\ast }\right\rangle \right\vert \\
&\leq &C(\int_{B_{X^{\#}}}\left\vert f\left( x\right) -f\left( x^{\prime
}\right) \right\vert ^{r}d\mu )^{\frac{1}{r}}\left\Vert y^{\ast }\right\Vert
_{L_{s}\left( \nu \right) } \\
&\leq &C(\int_{B_{X^{\#}}}\left\vert \left( \tilde{x}-\tilde{x}^{\prime
}\right) \left( f\right) \right\vert ^{r}d\mu )^{\frac{1}{r}}\left\Vert
y^{\ast }\right\Vert _{L_{s}\left( \nu \right) } \\
&\leq &C\left\Vert \tilde{x}-\tilde{x}^{\prime }\right\Vert _{L_{r}\left(
\mu \right) }\left\Vert y^{\ast }\right\Vert _{L_{s}\left( \nu \right) },
\end{eqnarray*}%
therefore by Theorem 2.3, $\overline{T}$ is Lipschitz-Cohen strongly $%
s^{\ast }$-summing.

\noindent $\left( 3\right) \implies \left( 2\right) :$ Given a factorization 
$\left( 3.3\right) ,$ consider the map $\delta _{\tilde{X}}\circ
T_{1}:X\longrightarrow \mathcal{F}\left( \tilde{X}\right) $ and $\widehat{%
T_{2}}:\mathcal{F}\left( \tilde{X}\right) \longrightarrow Y\subset Y^{\ast
\ast },$ then $T=\widehat{T_{2}}\circ \delta _{\tilde{X}}\circ T_{1},$ so,
we obtain what we needed for $\left[ \text{2, Theorem 5.4 (c)}\right] .\quad
\blacksquare 
\vspace{0.5cm}%
$

\textbf{Acknowledgments}

First, I am very grateful to the referee of Journal of Functional Analysis
for several interesting remarks and comments in his report that I have
followed which improved the first version of this paper, I am very
recognized for him. I thank also the referee of Journal JMAA for several
valuable suggestions which improve the second and final version of this
paper (Thanks for all).


\begin{thebibliography}{99}
\bibitem{} \textsc{R. Arens, J. Eells,} \textit{On embedding uniform and
topological spaces}, Pacific J. Math. 6 (1956) 397-403.

\bibitem{} \textsc{J. A. Chavez Dominguez,} \textit{Duality for Lipschitz }$%
p $\textit{-summing operators}. Journal of functional Analysis\textit{.} 01
2011, 261(2) 387-407.

\bibitem{} \textsc{S. Chevet,} \textit{Sur certains produits tensoriels
topologiques d'espaces de Banach}, Z. Wahrscheinlichkeitstheorie verw. Geb.
11, (1969) 120-138.

\bibitem{} \textsc{J. S. Cohen,} \textit{Absolutely }$p$-\textit{summing, }$%
p $-\textit{nuclear operators and their conjugates,} Math. Ann. 201 (1973)
177-200.

\bibitem{} \textsc{J. Diestel, H. Jarchow and A. Tonge,} \textit{Absolutely
summing operators,} Cambridge University Press, 1995.

\bibitem{} \textsc{M. Dubei, E.D. Tymchatyn and A. Zagorodnyuk,} \textit{%
Free Banach spaces and extension of Lipschitz maps}. Topology 48, (2009)
203-213.

\bibitem{} \textsc{J. D. Farmer and W. B. Johnson,} \textit{Lipschitz
p-summing Operators}, Proc. Amer. Math. Soc. 137, no.9, (2009) 2989-2995.

\bibitem{} \textsc{J. Flood,} Free Topological Vector Spaces, Dissertationes
Math. (Rozprawy Mat.) 221, 1984.

\bibitem{} \textsc{G. Godefroy and N. J. Kalton,} \textit{Lipschitz-free
Banach spaces}, Studia Mathematica 159 (1) (2003) 121-141.

\bibitem{} \textsc{A. Grothendieck,} R\'{e}sum\'{e} de la th\'{e}orie m\'{e}%
trique des produits tensoriels topologiques, Boll. Soc. Mat. S\~{a}o-Paulo 8
(1953) 1-79.

\bibitem{} \textsc{A. Jimenez-Vargas, J. M. Sepulcre, and Moises
Villegas-Vallecillos}, \textit{Lipschitz compact operators}, J. Math. Anal.
Appl. 415 (2014), no. 2, 889-901.

\bibitem{} \textsc{V.G. Pestov,} \textit{Free Banach spaces and
representation of topological groups}, Funct. Anal. Appl. 20 (1986) 70-72.

\bibitem{} \textsc{A. Pietsch,} Operators ideals, North-Holland, 1980.

\bibitem{} \textsc{A. Pietsch,} Absolut $p$-summierende Abbildungen in
normierten R\"{a}umen, Studia Math. 28 (1967) 333-353.

\bibitem{} \textsc{P. Saphar,} \textit{Applications \`{a} puissance nucl\'{e}%
aire et applications de Hilbert-Schmidt dans les espaces de Banach}, Annales
scientifiques de L'E.N.S. $3^{e}$ s\'{e}rie, tome 83, $n%
{{}^\circ}%
$ $2$ (1966) 113-151.

\bibitem{} \textsc{N. Weaver,} Lipschitz Algebras, World Scientific,
Singapore, New Jersey, London, Hong Kong, 1999.
\end{thebibliography}
\end{document}